\documentclass[a4paper,10pt]{article}
\usepackage{paper-en}
\usepackage{hyperref}

\def\thetitle{(Sym)metric planespotting}
\def\theauthors{Nina Lebedeva and Anton Petrunin}

\hypersetup{colorlinks=true,
citecolor=black,
linkcolor=black,
anchorcolor=black,
filecolor=black,
menucolor=black,
urlcolor=black,
pdftitle={\thetitle},
pdfauthor={\theauthors}
}

\begin{document}

\title{\thetitle}
\author{\theauthors}
\date{}
\maketitle

\begin{abstract}
We present a way to characterize the Euclidean and Lobachevsky planes metrically, starting with geodesic segments and imposing conditions on isometries.
\end{abstract}

We assume that the reader is not afraid of metric spaces.

\begin{thm}{Theorem}\label{thm:main}
A nontrivial 3-point-homogeneous geodesic space in which any two distinct points are fixed by exactly one nonidentity isometry is isometric to a Euclidean or Lobachevsky plane.
\end{thm}

We allow a Lobachevsky plane to have any constant negative curvature, so rescaling one gives another.
Let us spell out the other terms.

A space is \emph{nontrivial} if it has at least two points.

A \emph{geodesic} in a metric space $M$ is a distance-preserving map from a real interval to $M$.
The image of a geodesic with endpoints $x$ and $y$ will be denoted by $[xy]$; we also call it a \emph{geodesic}.
If any two points in $M$ can be joined by a geodesic, then $M$ is called a \emph{geodesic space}.

An \emph{isometry} of $M$ is a bijection from $M$ to itself that preserves distances.
Three-point homogeneity is essentially the side-side-side congruence criterion.
A metric space $M$ is \emph{3-point-homogeneous} if for any points $x_1$, $x_2$, $x_3$, $y_1$, $y_2$, $y_3\in M$ the equalities $|x_i-x_j|=|y_i-y_j|$ for all $i$ and $j$ imply that some isometry of $M$ sends $x_i$ to $y_i$ for each~$i$.
Here $|p-q|$ denotes the distance between points $p,q\in M$.
The points in a triple need not be distinct; in particular, 3-point homogeneity implies 2-point homogeneity and ordinary (1-point) homogeneity.

In a Euclidean or Lobachevsky plane, the only nonidentity isometry fixing two distinct points is reflection in the line thru them.
Thus the last condition of the theorem can be regarded as the \emph{unique reflection axiom}.
The round sphere is also 3-point-homogeneous, but its antipodal points are fixed by infinitely many isometries.

Since doubling distances changes the curvature of a Lobachevsky plane, adding the following \emph{self-similarity axiom} produces a characterization of the Euclidean plane:

\smallskip

\noindent\hskip2.5em\textit{The space is isometric to a copy of itself with all distances doubled.}

\medskip

The theorem provides a metric axiomatization of plane geometry.
A similar system was developed by Birger Iversen \cite{iversen}.
Both systems are built on metric spaces and their symmetries.
Iversen assumes that any two distinct points lie on a unique line (an isometric copy of $\RR$), and we prove this --- this is the main difference.

Our proof is based on Iversen's classification theorem, which is stated at the very end.
This note was prompted by another result of ours~\cite{lebedeva-petrunin}.

\section*{Proof}

Assume that a metric space $M$ satisfies the assumptions of our theorem.

{

\begin{wrapfigure}{r}{35mm}
\centering
\vskip-2mm
\includegraphics{mppics/pic-10}
\vskip-0mm
\end{wrapfigure}

\parbf{Two distances.}
Let $p,q\in M$ be two distinct points.
The distances $r=|p-x|$ and $s=|q-x|$ leave at most two possibilities for $x$.

Indeed, there is a unique nonidentity isometry, say~$\sigma$, fixing $p$ and~$q$.
If $|p-y|=r$ and $|q-y|=s$, then 3-point homogeneity gives an isometry fixing $p$ and $q$ and sending $x$ to $y$.
Thus $y=x$ if the isometry is the identity, and $y=\sigma(x)$ otherwise --- there is no third option.

}

\parbf{Local uniqueness of geodesics.}
We will show that there is $\eps_0>0$ such that any two points $x$ and $y$ at distance less than $\eps_0$ are joined by a unique geodesic.

A \emph{midpoint} of $x$ and $y$ is a point $m$ such that
\[|x-m|=|y-m|=\tfrac12\cdot|x-y|.\]
If the midpoint is unique, then every geodesic $[xy]$ contains it.
If midpoints are unique for all pairs of points at distance less than $\eps_0$, then repeated bisection gives uniqueness of geodesics of length less than $\eps_0$.

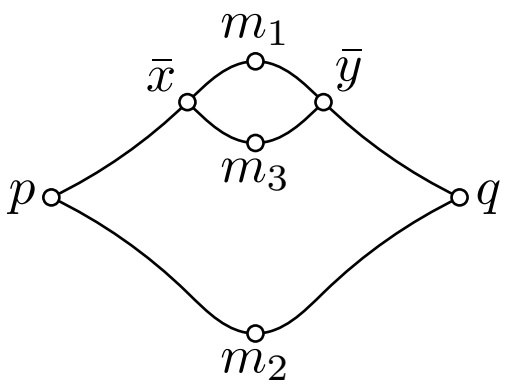
{

\begin{wrapfigure}{r}{35mm}
\centering
\vskip-5mm
\includegraphics{mppics/pic-20}
\vskip-0mm
\end{wrapfigure}

As observed above, a pair of points has at most two midpoints.
Suppose midpoints are not always unique, say $p$ and $q$ have two midpoints $m_1$ and $m_2$; set
\[\eps_0\z=|m_1-m_2|.\]

Let $0<|x-y|<\eps_0$.
Choose points $\bar x\z\in [m_1p]$ and $\bar y\in [m_1q]$ at distance $\tfrac12\cdot|x-y|$ from $m_1$.
Since $\bar x$ and $\bar y$ lie on a geodesic from $p$ to $q$ passing thru $m_1$, we have $|\bar x-\bar y|=|x-y|$.

}

Suppose $x$ and $y$ have two midpoints, say $n_1$ and $n_2$.
By 3-point homogeneity, there is an isometry
$\iota$ sending $x\mapsto \bar x$, $y\mapsto \bar y$, and $n_1\mapsto m_1$.
By the triangle inequality, $m_3=\iota(n_2)$ is also a midpoint of $p$ and $q$.
Note that $m_3\ne m_1$ and $m_3\ne m_2$; the latter follows since
\[|m_3-m_1|\le |x-y|<\eps_0=|m_1-m_2|.\]
Thus the midpoints $m_1$, $m_2$, and $m_3$ are distinct, which is impossible.

\parbf{No branching.}
We will show that \textit{geodesics in $M$ do not branch}, that is, if
two geodesics $[xz]$ and $[xz']$ contain a nontrivial geodesic $[xy]$, then $[xz]$ contains $[xz']$ or $[xz]$ is contained in $[xz']$.

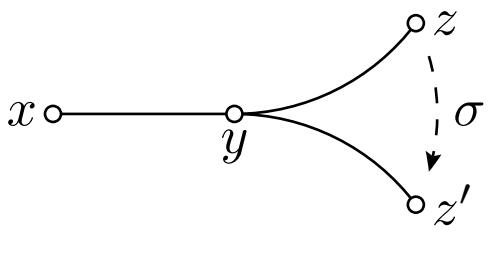
{

\begin{wrapfigure}{r}{35mm}
\centering
\vskip-4mm
\includegraphics{mppics/pic-30}
\vskip-0mm
\end{wrapfigure}

If geodesics branched, we could shorten them and move $x$ closer to the branching point so that $|x-z|\z=|x-z'|\z<\eps_0$ and $z\ne z'$.

By 3-point homogeneity, there is an isometry $\sigma$ fixing $x,y$ and sending $z$ to~$z'$.
But there is also a nonidentity isometry, say $\tau$, fixing $x$ and $z$; it must fix every point on the unique geodesic $[xz]$, including~$y$.
Thus $\sigma$ and $\tau$ are two distinct nonidentity isometries fixing $x$ and $y$ --- a contradiction.

}

\parbf{Lines.}
Any geodesic can be extended to a \emph{local geodesic} $\gamma\colon\mathbb R\to M$; this means that sufficiently short segments of $\gamma$ are geodesics.
For this, use 2-point homogeneity to align a short final segment with the start of a longer geodesic.
By local uniqueness, they agree on the overlap, giving an extension.
Repeating this with the same small step in both directions produces $\gamma$.
The absence of branching ensures that a local geodesic is uniquely determined by any of its nontrivial segments.

Let $d$ be the diameter of $M$; we will show that $d=\infty$.
The following argument is borrowed from Garrett Birkhoff~\cite{birkhoff1944}.

Note that $M$ contains a geodesic of any length less than $d$.
By 2-point homogeneity and the absence of branching, these geodesics can be aligned with the local geodesic $\gamma$ constructed above.
Thus
\[|\gamma(t)-\gamma(s)|=|t-s|\quad\text{if}\quad |t-s|<d.
\eqlbl{eq:|t-s|}\]

Suppose $d<\infty$; then by continuity, \ref{eq:|t-s|} also holds when $|t-s|=d$.
Choose a small $\eps>0$ and set
\[p=\gamma(0),
\quad x=\gamma(d+\eps),
\quad r=|p-x|,
\quad y=\gamma(r),
\quad\text{and}\quad
z=\gamma(\eps).\]
By the triangle inequality, $d-\eps\le r\le d$.

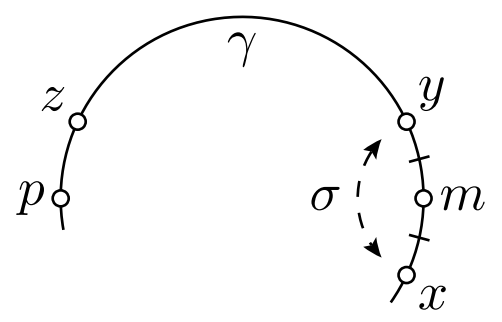
\begin{wrapfigure}{r}{35mm}
\centering
\vskip-8mm
\includegraphics{mppics/pic-40}
\vskip-0mm
\end{wrapfigure}

The points $x$ and $y$ are equidistant from $p$, so there is an isometry $\sigma$ fixing $p$ and swapping $x$ and~$y$.
It reverses the unique geodesic $[xy]$ and hence fixes its midpoint, say $m$.
By our assumptions, there is another nonidentity isometry fixing $p$ and $z$;
it must fix the unique geodesic $[pz]$, and hence every point on $\gamma$.
We now have two distinct nonidentity isometries fixing $p$ and $m$ --- a contradiction.

Thus $d=\infty$, and every segment of the local geodesic $\gamma$ is a geodesic;
that is, $\gamma$ is a geodesic, a distance-preserving embedding of the whole real line into $M$.
The image of such a geodesic will be called a \emph{line}.

Let us check that the geodesic between points $x$ and $y$ is unique.
Suppose there are two.
Extend one of them beyond $y$ to a point $z$.
The other, together with the same segment $[yz]$, is also a geodesic $[xz]$.
The two geodesics branch, which is impossible.

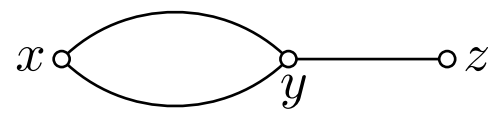
\begin{wrapfigure}{r}{35mm}
\centering
\vskip-5mm
\includegraphics{mppics/pic-50}
\vskip-0mm
\end{wrapfigure}

Thus there is a unique line thru any two distinct points.
There is also a unique nonidentity isometry fixing every point on a given line; we call it the \emph{reflection} in that line.
The assumptions of the theorem also imply that a reflection is an involution, meaning that it is its own inverse.

\parbf{Half-planes.}
The first part of the following argument was already used by Benjamin Kagan~\cite[Satz~30]{kagan1902}.

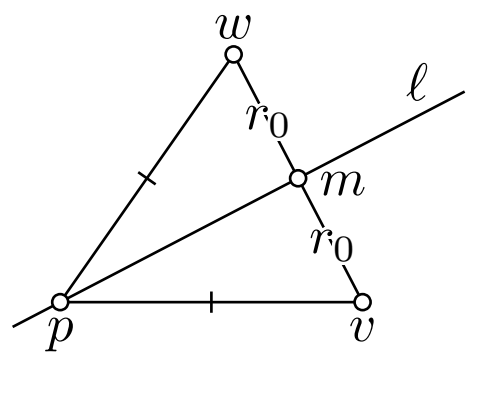
\begin{wrapfigure}[6]{r}{35mm}
\centering
\vskip-6mm
\includegraphics{mppics/pic-120}
\vskip-0mm
\end{wrapfigure}

We can choose two geodesics $[pv]$ and $[pw]$ of equal length on distinct lines thru $p$.
Indeed, since $M$ is nontrivial, it contains a line.
The reflection in this line moves some point, which must lie off the line.

Let $m$ be the midpoint of $[vw]$.
By 3-point homogeneity, there is an isometry $\sigma$ fixing $p$ and swapping $v$ and~$w$.
It fixes $m$ and hence fixes the line $pm$ pointwise; denote this line by $\ell$.

We claim that $\ell$ consists of all points equidistant from $v$ and $w$.
On each half-line of $\ell$ from $m$, the common distance to $v$ and $w$ starts at $r_0=\tfrac12\cdot|v-w|$ and tends to infinity.
Thus for each $r>r_0$ there is a point on each half-line at distance $r$ from both $v$ and $w$.
Since two distances to $v$ and $w$ leave only two possibilities, there is no third.
For $r=r_0$ there is only the midpoint $m$, and $r<r_0$ is ruled out by the triangle inequality.

Let
\[V=\set{x\in M}{|x-v|<|x-w|}.\]
If $x\in V$, then $[vx]\subset V$ by the triangle inequality.
Thus $V$ and $\sigma(V)$ are connected disjoint open sets covering $M\setminus\ell$, so they are its two components, and $\sigma$ swaps them.

By 2-point homogeneity, an isometry of $M$ can move $\ell$ to any given line.
Thus $M$ meets the assumptions in the following theorem of Birger Iversen~\cite[p.~ix]{iversen}.

\begin{thm}{Classification theorem}
A nontrivial metric space $M$ is isometric to a Euclidean or Lobachevsky plane if it meets the following conditions:

\smallskip

\noindent\textsc{incidence axiom.}
Thru two distinct points of $M$ there passes a unique line.

\smallskip

\noindent\textsc{reflection axiom.}
The complement of a given line in $M$ has two connected
components, and there exists an isometry $\sigma$ of $M$ which fixes the points of the line, but interchanges the two connected components of its complement.

\end{thm}

\parit{We used AI (GPT-6.2 Astra) to discuss the proof, search the literature, translate, and edit the text.}

{\sloppy
\def\emph{\textit}
\printbibliography[heading=bibintoc]\fussy}
\end{document}